%% file: turnpike.tex
\documentclass{amsart}
\pdfoutput=1 

\usepackage{amssymb}
\usepackage{amsthm}

\usepackage{tikz}
\usepackage{verbatim}
\usepackage{pgfplots}
\usetikzlibrary{arrows,positioning,decorations}
\tikzset{
    >=stealth',
    punkt/.style={
           rectangle,
           rounded corners,
           draw=black, very thick,
           text width=6.5em,
           minimum height=2em,
           text centered},
    pil/.style={
           ->,
           shorten <=2pt,
           shorten >=2pt,}
}

\makeatletter
\@namedef{subjclassname@2020}{%
  \textup{2020} Mathematics Subject Classification}
\makeatother

\newtheorem{theorem}{Theorem}
\newtheorem{lemma}[theorem]{Lemma}
\newtheorem{corollary}[theorem]{Corollary}
\newtheorem{proposition}[theorem]{Proposition}

\newtheorem{definition}[theorem]{Definition}

\theoremstyle{definition}
\newtheorem{example}[theorem]{Example}

\newtheorem*{acknowledgment}{Acknowledgment}

\newcommand{\bR}{\mathbb{R}}

\newcommand{\bN}{\mathbb{N}}
\newcommand{\bZ}{\mathbb{Z}}

\newcommand{\cC}{\mathcal{C}}
\newcommand{\C}{\cC}

\newcommand{\cI}{\mathcal{I}}
\newcommand{\I}{\cI}
\newcommand{\cJ}{\mathcal{J}}
\newcommand{\J}{\cJ}

\newcommand{\fin}{\textrm{Fin}} 

\newcommand{\card}[1]{{\#{}#1}} 
\newcommand{\norm}[1]{\left\|#1\right\|} 

\newcommand{\imply}{\Rightarrow}

\newcommand{\probA}{(*)}
\newcommand{\probB}[1]{(#1{\rm/}\!**)}
\newcommand{\condB}[1]{(#1{\rm/}C2)}

\begin{document}

\title{Invariant ideals and its applications to the turnpike theory}

\author[M. Mammadov]{Musa Mammadov} 
\address{School of Information Technology, Deakin University, VIC 3125, Australia}
\email{musa.mammadov@deakin.edu.au}
\author[P. Szuca]{Piotr Szuca}
\address{Institute of Mathematics, University of Gda\'{n}sk, ul.~Wita
  Stwosza 57, 80-952 Gda\'{n}sk, Poland}
\email{Piotr.Szuca@ug.edu.pl}

\date{\today}

\begin{abstract}
In this paper the turnpike property is established for a non-convex optimal control problem in discrete time. 
The functional is defined by the notion of the ideal convergence and can be considered as an analogue of the terminal functional defined over infinite time horizon.
The turnpike property states that every optimal solution converges to some unique optimal stationary point
in the sense of ideal convergence 
if the ideal is invariant under translations.
This kind of convergence generalizes, for example, statistical convergence and convergence with respect to logarithmic density zero sets.
\end{abstract}

\subjclass[2020]{Primary 40A35; Secondary 49J99, 54A20}

\keywords{$\I$-convergence, $\I$-cluster set, statistical convergence,
  turnpike property, optimal control, discrete systems}


\maketitle


\input{intro.tex}

\input{ideal.tex}

\input{problem.tex}

\input{main.tex}

\input{proofs.tex}


\begin{acknowledgment}
The authors are indebted to Paolo Leonetti for his critical reading of the manuscript.
\end{acknowledgment}

\bibliographystyle{plain}
\bibliography{turnpike}

\end{document}

%% file: intro.tex
\section{Introduction}

The turnpike theory investigates an important property of dynamical systems. It can be considered as a theory
that justifies the importance of some equilibrium/stationary states.
For example, in macroeconomic models
the turnpike property states that, \emph{regardless of initial
conditions, all optimal trajectories spend most of the time
within a small neighborhood of some optimal stationary point when the planning period is long enough.}
Obviously, in the absence of such a property, using some of optimal stationary points
as a criteria for ``good'' policy formulation might be misleading.
Correspondingly, the turnpike property is in the core of many important
theories in economics.

Many real-life processes are happening in an optimal way and have the tendency to stabilize;
that is, the turnpike property is expected to hold for a broad class of problems.
It provides valuable insights into the nature of these processes
by investigating underlying principles of evolution that lead to stability.
It can also be used to assess the ``quality'' of mathematical modeling and to develop more adequate
equations describing system dynamics as well as optimality criteria.

The first result in this area is obtained by John von Neumann (\cite{Neu45}) for discrete
time systems. The phenomenon is called the turnpike property after Chapter 12, \cite{dss1958} by
Dorfman, Samuelson and Solow. For a classification of different definitions for this property, see \cite{CHL91,mr1977,McK76,Zas06},
as well as \cite{dgsw2012} for the so-called \emph{exponential} turnpike property.
Possible applications in Markov Games can be found in a recent study \cite{Kolo12}.

The approaches suggested for the study of the turnpike property involve continuous and discrete time systems.
Some convexity assumptions are sufficient for discrete time systems
\cite{mr1977,McK76}; however, rather restrictive assumptions are
usually required for continuous time systems.
The majority of them deal with the (discounted and undiscounted) integral functionals.
We mention here the approaches developed by Rockafellar \cite{Roc73,Roc76}, Scheinkman, Brock and collaborators (see, for example, \cite{Sch79, Sch1976}),  Cass and Shell \cite{cass1976}, Leizarowitz \cite{Leiz89}, Mamedov \cite{Mam-1992},  Montrucchio \cite{montrucchio1995turnpike},
Zaslavski \cite{Zas11,Zas132,Zas13} (we refer to \cite{CHL91,Zas06} for more references).

In this paper we consider an optimal control problem in discrete time. It
extends the results obtained in \cite{Mam-1985} where
a special class of terminal functionals is introduced as a lower limit at infinity of
utility functions. This approach allowed to establish the turnpike property for
a much broader class of optimal control problems than those involving integral
functionals (discounted and undiscounted).

Later, this class of
terminal functionals was used to establish a connection between the
turnpike theory and the notion of statistical convergence
\cite{MamPeh-2001-JMAA, PehMam-2000-Opt}; as a result, the convergence of optimal
trajectories is proved in terms of the statistical
(``almost'') convergence.
These terminal functionals also allowed the extension of the turnpike theory to time delay systems;
the first results in this area have been established in several
recent papers \cite{IvaMamTr-2012-JOGO,Mam-2013-SIAM}.
Moreover, some generalizations based on the notion of the $A$-statistical cluster points have been
obtained in \cite{MR3166601}.

The main purpose of this paper is to formulate the optimality criteria
by using the notion of ideal convergence. As detailed in the next section, the ideal convergence is a more general concept than
the statistical convergence as well as the $A$-statistical convergence.
In this way the turnpike property is established for a broad class of
non-convex optimal control problems where the asymptotical stability of optimal trajectories is
formulated in terms of
the ideal convergence.

Recently (and independently) Leonetti and Caprio in~\cite{MR4428911} considered turnpike property for ideals invariant under translation
in the context of normed vector spaces. We discuss our approaches in Section~\ref{sec:main}.

The rest of this article is organized as follows. In the next section
the definition of the ideal, its properties and some particular cases, including the statistical convergence, are provided.
In Section \ref{sec:problem} we formulate the optimal control problem and main assumptions.
The main results of the paper --- the turnpike theorems are provided in Section \ref{sec:main}.
The proof of the main theorem is in Section \ref{sec:proofs}.

%% file: ideal.tex
\section{Convergence with respect to ideal vs statistical convergence}


Let $x=(x_n)_{n\in\bN}$ be a sequence of elements of $\bR^m$.
For the sake of simplicity, we will consider the Euclidean norm $\norm{\cdot}.$
The classical definition of convergence of $x$
to $a$ says that for every $\varepsilon>0$ the set of all $n\in\bN$
with $\norm{x_n-a}\geq\varepsilon$ is finite, i.e. it is ``small'' in some sense.
If we understand the word ``small'' as ``of asymptotic density zero'' then we obtain
the definition of statistical convergence (Def.~\ref{def:statistical-convergence}). The same method can be used
to formulate the definition of statistical cluster point.
The classical one says that $a$ is a cluster point of $x$ if
for every $\varepsilon>0$ the set of all $n\in\bN$
with $\norm{x_n-a}<\varepsilon$ is infinite, i.e. it has ``many'' elements.
If ``many'' means
``not of asymptotic density zero'' then we obtain
the definition of statistical cluster point (see e.g.~\cite{fridy}).

One of the possible generalizations of this kind of being ``small'' (having ``many'' elements) is ``belonging to the ideal''
(``be an element of co-ideal'').

The cardinality of a set $X$ is denoted by $\card{X}$.
$\mathcal{P}(\bN)$ denotes the power set of $\bN$.

\begin{definition}\label{def:ideal-convergence}
An \emph{ideal on $\mathcal{P}(\bN)$}
is a family $\I\subset\mathcal{P}(\bN)$
which is non-empty, hereditary and closed under taking finite unions,
i.e. it fulfills the following three conditions:
\begin{enumerate}
\item $\emptyset\in\I$;
\item $A\in\I$ if $A\subset B$ and $B\in\I$;
\item $A\cup B\in\I$ if $A,B\in\I$.
\end{enumerate}
\end{definition}

\begin{figure}\centering
\begin{tikzpicture}[node distance=1cm, auto,]
 \node (fin) {$\fin$};
 \node[right=1cm of fin] (dummy) {};
 \node[above=0.5cm of dummy] (w) {$\mathcal{W}$} edge[pil,<-] (fin);
 \node[below=0.5cm of dummy] (sum) {$\I_{1/n}$} edge[pil,<-] (fin);
 \node[right=1cm of dummy] (density) {$\I_d$} edge[pil,<-] (w) edge[pil,<-] (sum);
 \node[right=1cm of density] (log) {$\I_{log}$} edge[pil,<-] (density);
\end{tikzpicture}
\caption{Inclusions of ideals, implications between $\I$-convergence,
  and inclusions of sets of $\I$-cluster points for ideals from Example~\ref{ex:ideals}.
  Arrow ``$\I\longrightarrow\J$'' means that ``$\I\subset\J$'', and for every sequence $x$,
  ``$x\to_\I a \imply x\to_\J a$'',
  ``$\Gamma_\I(x)\supset\Gamma_\J(x)$''.}\label{fig:ideal-implications}
\end{figure}
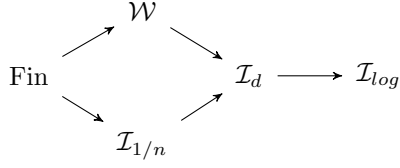

\begin{example}\label{ex:ideals}
By $\fin$ we denote the ideal of all finite
subsets of $\bN=\{1,2,\ldots\}$.
There are many examples of ideals considered in the literature, e.g.
\begin{enumerate}
\item[$(1)$]
the ideal of sets of asymptotic density zero
$$\I_d=\left\{A\subset\bN\ :\ \overline{d}(A)=0\right\},$$
where $\overline{d}\colon\mathcal{P}(\bN)\to[0,1]$ is given by the formula
$$\overline{d}(A)=\limsup_{n\to\infty} \frac{\card{(A\cap\{1,2,\ldots,n\})}}{n}$$
is the well-known definition of upper asymptotic density of the set $A$;
\item[$(2)$]
the ideal of sets of logarithmic density zero
$$\I_{log}=\left\{A\subset\bN\ :\
    \limsup_{n\to\infty}\frac{\sum_{k\in A\cap\{1,2,\ldots,n\}}\frac{1}{k}}{\sum_{k\leq n}\frac{1}{k}}=0\right\};$$
\item[$(3)$]
the ideal
$$\I_{1/n}= \left\{A\subset\bN: \sum_{n\in
    A}\frac{1}{n}<\infty\right\};$$
\item[$(4)$]
the ideal of arithmetic progressions free sets
$$\mathcal{W}=\left\{W\subset\bN\ :\ \textup{$W$ does not contain
      arithmetic progressions of all lengths}\right\}.$$
\end{enumerate}
\end{example}
Ideals $\I_d$ and $\I_{log}$ belongs to the wider class of Erd\H{o}s-Ulam ideal's
(defined by submeasures of special kind, see~\cite{just-krawczyk}).
Ideal $\I_{1/n}$ is an representant of the class of summable ideals (see~\cite{mazur}).
The fact that $\mathcal{W}$
is an ideal follows from the non-trivial theorem of van der Waerden
(this ideal was considered by Kojman
in~\cite{waerden-spaces}).
One can also consider trivial ideals $\I=\mathcal{P}(\bN)$,
$\I=\{\emptyset\}$,
or principal ideals $\I_n=\{A\subset\bN\colon n\notin A\}$,
however they are not interesting from our point of view.
If not explicitly said
we assume that
all considered ideals are proper (i.e.~$\I\not=\mathcal{P}(\bN)$) and contain all finite sets
(i.e.~$\fin\subset\I$).
The inclusions between abovementioned families are shown on Figure~\ref{fig:ideal-implications}.
The only non-trivial inclusions are: $\I_{1/n}\subset\I_d$ (a folklore application of Cauchy condensation test),
$\mathcal{W}\subset\I_d$ (the famous theorem of Szemer\'{e}di),
and $\I_d\subset\I_{log}$ (by well-known
inequalities between upper logarithmic density and upper asymptotic density).
It is easy to observe that $\I_{1/n}\not\subset\mathcal{W}$,
but the status of the inclusion $\mathcal{W}\subset\I_{1/n}$
is unknown (``Erd\H{o}s conjecture
  on arithmetic progressions''
  says that the van der Waerden ideal $\mathcal{W}$ is contained in the
  ideal $\I_{1/n}$.)


\subsection{$\I$-convergence and $\I$-cluster points}

The notion
  of the ideal convergence is dual (equivalent) to the
  notion of the filter convergence introduced by Cartan in 1937
  (\cite{cartan}). The notion of the filter convergence
  has been an important tool in general topology and functional
  analysis since 1940 (when Bourbaki's book~\cite{Bourbaki} appeared).
  Nowadays many authors prefer to use an equivalent dual notion of the
  ideal convergence (see e.g.~frequently quoted work \cite{ksw}).

\begin{definition}
A sequence $(x_n)_{n\in\bN}$ of elements of $\bR^m$ is said to be
\emph{$\I$-convergent} to $a\in\bR^m$ ($a = \I-\lim x_n$, or $x_n\to_\I a$, in short) if and only if for
each $\varepsilon>0$
$$\left\{n\in\bN\ :\ \norm{x_n-a}\geq\varepsilon\right\}\in\I.$$
\end{definition}

The sequence $(x_n)$ is convergent to $a$ if and only if it is $\fin$-convergent to $a$.

It is also easy to see that for any sequence $x=(x_n)$ and two ideals $\I$, $\J$, if $\I\subset\J$ then
$x\to_{\I} a$ implies that $x\to_{\J} a$ (see Figure~\ref{fig:ideal-implications}).

\begin{definition}
The $a\in\bR^m$ is an \emph{$\I$-cluster point} of a sequence $x=(x_n)_{n\in\bN}$ of elements of $\bR^m$
if for each $\varepsilon>0$
$$\{n\in\bN\colon \norm{x_n-a}<\varepsilon\}\notin\I.$$
\end{definition}

By \emph{$\I$-cluster set of $x$} we understand the set
$$\Gamma_\I(x)=\left\{a\in\bR^m\colon a\textup{\ is an\ }\I\textup{-cluster point of\ }x\right\}.$$
Recall that $\Gamma(x)=\Gamma_{\fin}(x)$ is a set of classical cluster (limit) points of $x$.

\begin{proposition}\label{prop:properties-of-I-convergence}
For any bounded sequence $x=(x_n)$,
\begin{enumerate}
\item[$(1)$] $\Gamma_\I(x)\not=\emptyset$ (\cite{generalized-stat-convergence}), and
\item[$(2)$] $\Gamma_\I(x)$ is closed (\cite{ksw}), and
\item[$(3)$] $\Gamma_\I(x)=\{a\}$ if and only if $x\to_\I a$.
\end{enumerate}
Moreover, if $\I\subset\J$ then
$\Gamma_\J(x)\subset\Gamma_\I(x)$
(\cite{generalized-stat-convergence}, see Figure~\ref{fig:ideal-implications}).
\end{proposition}

Part~$(3)$ follows from the folklore argument:
$a$ is the unique $\I$-cluster point of $x$, iff
$\{n\colon ||x_n-a||\geq\varepsilon\}\in\I$ for every $\varepsilon>0$, iff
$x_n\to_\I a$.


\subsection{$\I$-convergence vs statistical convergence}

The notion of the ideal convergence is a common
generalization of the classical notion of convergence and statistical convergence.
The concept of statistical convergence was introduced by Fast \cite{fast} and then
it was studied by many authors.

\begin{definition}[\cite{fast}]\label{def:statistical-convergence}
A sequence $x=(x_n)_{n\in\bN}$ of elements of $\bR^m$ is said to be \emph{statistically convergent
to an $a\in\bR^m$} if for each $\varepsilon>0$
the set of all indices $n$ such that
$\{n\in\bN\colon \norm{x_n-a}\geq\varepsilon\}$ has upper asymptotic density zero, i.e.
$$\overline{d}(\{n\in\bN\colon \norm{x_n-a}\geq\varepsilon\})=0, \textup{\ for all\ } \varepsilon>0.$$
\end{definition}

Obviously,
$x$ is statistically convergent to $a$ if and only if $x\to_{\I_d}a$.
Following the concept of a statistically convergent sequence Fridy in~\cite{fridy}
introduced the notion of a statistical cluster point, which---using our notation---is equal
to the notion of $\I_d$-cluster point.
%
%
Proposition~\ref{prop:properties-of-I-convergence} in case of statistical convergence
was proved in \cite{fridy}.
Since
statistical convergence is a particular case of $\I$-convergence,
each theorem which has an ideal variant is also true
in its statistical version. However, in the sequel
we will use some lemmas which were formulated in the literature
for the case of statistical convergence and statistical cluster points.

An open $\varepsilon$-neighbourhood of a given set $A\subset\bR^m$
will be denoted by
$$ B(A,\varepsilon) = \{y\in\bR^m\colon \exists_{a\in A} \norm{a-y}<\varepsilon\}.$$
For each $a\in\bR^m$ we do not distinguish between $B(\{a\},\varepsilon)$ and $B(a,\varepsilon)$.

\begin{lemma}[\cite{PehMam-2000-Opt}]\label{le:st-cluster}
Let $x = (x_k)_{k\in\bN}$ be a bounded sequence.
Then for any $\varepsilon > 0$
$$\overline{d}(\left\{k \in\bN \colon x_k\notin B(\Gamma_{\I_d} (x),\varepsilon) \right\})=0.$$
\end{lemma}

The ideal version of the above lemma can be proved using the same method as in~\cite{PehMam-2000-Opt},
but we give a short proof using~\cite[Le.~3.1]{csst}.

\begin{lemma}[\cite{csst}]\label{le:cluster-compact}
Suppose that $\I$ is an ideal, $(x_n)$ is a sequence and $K\subset\bR^m$ is compact.
If $\{n\in\bN\colon x_n\in K\}\notin\I$ then $K\cap\Gamma_\I(x)\not=\emptyset$.
\end{lemma}

\begin{lemma}[Ideal version of Lemma~\ref{le:st-cluster}]\label{le:cluster-nbhd}
Let $x = (x_k)_{k\in\bN}$ be a bounded sequence.
Then for any ideal $\I$ and $\varepsilon > 0$
$$\left\{k \in\bN \colon x_k\notin B(\Gamma_{\I} (x),\varepsilon) \right\}\in\I.$$
\end{lemma}

{\bf Proof.~}
Since $x$ is bounded, there exists a compact set $C$ such that $x_n\in C$ for all $n$.
If we assume that $\left\{k \in\bN \colon x_k\notin B(\Gamma_{\I} (x),\varepsilon) \right\}\notin\I$,
then the set $K=C\setminus B(\Gamma_{\I} (x),\varepsilon)$ is compact
and $\{n\in\bN\colon x_n\in K\}\notin\I$. By Lemma~\ref{le:cluster-compact}, $K\cap\Gamma_\I(x)\not=\emptyset$,
a contradiction.
$\Box$ \bigskip


\subsection{Ideals invariant under translations}
\label{subsec:inv}

By $\bZ$ we denote the set of all integers.

\begin{definition}
We say that an ideal $\I$ is {\em invariant under translations}
if for each $A\in\I$ and $i\in\bZ$,
$$A+i\in\I,
\textup{\ where\ }
A+i=\left\{a+i\ :\ a\in A\right\}\cap\bN.$$
\end{definition}

All ideals considered in Example~\ref{ex:ideals} are invariant under translations.
For the proof of this fact and other examples see~\cite{filipow-szuca}.

Our main results from Section~\ref{sec:main} are valid for ideals invariant under translations.
The key argument for this fact is the following property of $\I$-cluster sets for such ideals.

\begin{lemma}\label{le:cluster-ideal-invariant}
Suppose that $\I$ is invariant under translations, $x=(x_k)_{k\in\bN}$ is a sequence in $\bR^m$.
Then for any non-empty $G\subset\Gamma_\I(x)$, $i\in\bZ$ and $\delta_1,\delta_2>0$:
$$\left\{k\in\bN\colon x_k\in B\left(G,\delta_1\right)\textup{\ and\ }x_{k+i}\in B\left(\Gamma_\I(x),\delta_2\right)\right\}\notin\I.$$
In particular, this set is non-empty.
\end{lemma}

{\bf Proof.~}
Let $K^1_{\delta_1}=\{k\in\bN\colon x_k\in B(G,\delta_1)\}$
and $K^2_{\delta_2}=\{k\in\bN\colon x_k\in B(\Gamma_\I(x),\delta_2)\}$.
Since $G\subset\Gamma_\I(x)$ and $G\not=\emptyset$, $K^1_{\delta_1}\notin\I$.
By Lemma~\ref{le:cluster-nbhd},
$\bN\setminus K^2_{\delta_2}\in\I$.
Consider the set
$K^1_{\delta_1}+i=\{k+i \colon k\in K^1_{\delta_1}\}$.
$\I$ is invariant under translations, so $K^1_{\delta_1}+i\notin\I$.
Let $K=(K^1_{\delta_1}+i) \cap K^2_{\delta_2}$.
Since $K$ is an intersection of two sets,
one from the coideal (i.e.~not from the ideal) and second from the dual filter (i.e.~its complement belongs to the ideal),
$K\notin\I$.
Consider the set $K-i=\{k-i\colon k\in K\}$.
Again, since $\I$ is invariant under translations, $K-i\notin\I$.
For each $k\in K-i$,
$x_k\in B\left(G,\delta_1\right)\textup{\ and\ }x_{k+i}\in B\left(\Gamma_\I(x),\delta_2\right)$.
$\Box$ \bigskip

If $\I$ is invariant under translations then either $\I$ is a trivial ideal $\{\emptyset\}$,
or $\I$ contains all finite sets (i.e.~$\fin\subset\I$).
Indeed, if there is a non-empty set $F\in\I$, then $\{n\}\in\I$ for each $n\in F$.
From the invariance of $\I$ it follows that $\{k\}\in\I$ for every $k\in\bN$.
Since $\I$ is closed on finite unions, each finite set belongs to $\I$.

%% file: problem.tex
\section{Optimal control problem and main assumptions}
\label{sec:problem}

Consider the problem
$$x_{n+1}=f(x_n,u_n), x_1=\zeta^0, u_n\in U, \eqno \probA$$
$$J_\I(x)=\I{\rm-}\liminf \phi(x_n) \rightarrow \max. \eqno \probB{\I}$$
Here $\zeta^0$ is a fixed initial point,
function $f\colon\bR^m\times\bR^t\to\bR^m$ is continuous, $U\subset\bR^t$ is a compact set,
$\phi\colon\bR^m\to\bR$ is a continuous function, and for any sequence of reals $y=(y_n)$
$$\I{\rm-}\liminf y=\sup\left\{y_0\in\bR\colon\left\{n\in\bN\colon y_n<y_0\right\}\in\I\right\}.$$

The pair $\langle u,x\rangle$ is called a process
if the sequences $x = (x_n)$ and $u = (u_n)$
satisfy $(*)$ for all $n\in\bN$
($x$ is called a trajectory and $u$ is called a control).

In the sequel we will use the following characterization of the functional $J_\I$
(\cite[Le.~4.1]{MR3166601}) that is a generalization of Lemma 3.1 in \cite{PehMam-2000-Opt} established for the statistical convergence,
as well as the corresponding result from \cite{L83} established for classical convergence (see also~\cite[Cor.~3.3]{MR4428911}).

\begin{lemma}\label{le:functional-j}
For any bounded trajectory $x = (x_n)_{n\in\bN}$ the following representation is true
$$J_\I(x)=\min\Gamma_\I(\phi(x))=\min_{\zeta\in\Gamma_\I(x)}\phi(\zeta).$$ 
\end{lemma}

We assume that there is a compact (bounded and closed)
set $C\subset\bR^m$ such that $x_n\in C$ for all trajectories;
that is, we assume that trajectories are uniformly bounded.

$\zeta\in\bR^m$ is called a stationary point if there exists $u_0\in U$ such that $f(\zeta,u_0) = \zeta$.
We denote the set of stationary points by $M$. It is clear that $M$ is a closed set.
$\zeta^\star\in M$ is called an optimal stationary point if 
$$\phi(\zeta^\star)=\phi^\star ~\dot{=}~ \max_{\zeta\in M}\phi(\zeta).$$
We will assume that the set of all optimal stationary points is non-empty.
This is not a restrictive assumption since function $\phi$ is continuous and
the set $M$ is closed; for example, it is satisfied if $M$ is in addition bounded.

Define the set
$$M^\star=\left\{\zeta^\star\in M\colon ~\zeta^\star \textup{\ is an optimal stationary point}\right\},$$
and
$$D^\star=\left\{\zeta\in C\colon ~\phi(\zeta)\geq\phi^\star\right\}.$$
We assume that the set $C$ is large enough to accommodate $M^\star;$ that is,
$M^\star \subset C$.
Then clearly, $M^\star=M\cap D^\star$.

Consider the following three conditions.
\begin{itemize}
\item[$(C1)$:]
  optimal stationary point $\zeta^\star$ is unique, i.e.~$M^\star=\{\zeta^\star\}$;
\item[$\condB{\I}$:]
  there exists a process $\langle u^\star,x^\star\rangle$ such that 
  $\Gamma_\I(x^\star)\subset D^\star$;
\item[$(C3)$:] there exists a continuous function $P\colon\bR^m\to\bR$
such that
$$
    P(f(x_0,u_0))<P(x_0)  ~~~~~ {\rm for~all}~~ x_0\in D^\star\setminus M^\star, u_0\in U,
$$
and
$$
    P(f(x_0,u_0))\leq P(x_0)  ~~~~~ {\rm for~all}~~ x_0\in D^\star, u_0\in U.
$$
\end{itemize}

One can also consider condition $(C2)=\condB{\fin}$:
\begin{itemize}
\item[$(C2)$:]
  there exists a process $\langle u^\star,x^\star\rangle$ such that any limit point of the sequence $x^\star$ is in $D^\star$.
\end{itemize}

Note that if the (unique) optimal stationary point $\zeta^\star$ belongs to the interior of $D^\star$, then proof of turnpike property
is not difficult and can be regarded as a ``trivial'' case where
condition $(C3)$ ensures the existence of some Lyapunov function,
with derivative $P,$ defined on a small neighborhood of  $\zeta^\star.$

The most interesting case is when an optimal stationary point $\zeta^\star$ belongs to a boundary of $D^\star$; that is, the
both sets $D^\star$ and $D^{\star-} = \left\{\zeta\in C\colon ~\phi(\zeta)<\phi^\star\right\}$
have nonempty intersection
with any small neighborhood of $\zeta^\star$.
In this case the inequality $P(f(x_0,u_0))>P(x_0)$ may hold for some
$x_0 \in D^{\star-}$; that is, condition $(C3)$ does not
guarantee the existence of Lyapunov functions.

Note also that condition $\condB{\I}/(C2)$ can be formulated equivalently  
``there exists a process $\langle u^\star,x^\star\rangle$ such that $J_\I(x^\star) \geq \phi^\star$'' (see \cite[$(A6)$]{MR4428911}),
or stronger ``there exists a process $\langle u^\star,x^\star\rangle$ such that $x^\star\to_\I \zeta^\star$'' 
(see \cite{Mam-1985, MamPeh-2001-JMAA, PehMam-2000-Opt}),

Recall that if $\I\subset\J$ then
$\I$-convergence is stronger than $\J$-convergence,
thus by 
Proposition~\ref{prop:properties-of-I-convergence}, $(C2)$ is stronger than $\condB{\I}$ for each non-trivial $\I$
which is invariant under translations.
Example~\ref{ex:statistical-not-fin} shows that these two conditions are really different,
i.e.~there exists a system for which $(C1)$, $\condB{\I_d}$, $(C3)$ hold, but
$(C2)$ does not hold (see also~\cite[Ex.~2.5]{MR4428911}).

\begin{example}\label{ex:statistical-not-fin}
Consider the middle-third Cantor set $T$.
It is homeomorphic to the space $\{0,1\}^\bN$
with the product (Tychonoff) topology;
for example, the formula
\begin{equation}\label{ex:def-of-homeomorphism}
\sum_{i=1}^\infty \frac{2\cdot a_i}{3^i} \textup{\ for any\ } a=\langle a_1,a_2,\ldots\rangle\in\{0,1\}^\bN
\end{equation}
gives us a homeomorphism between $\{0,1\}^\bN$ with Tychonoff topology and middle-third Cantor set.
In this example we will not distinguish between $T$ and $\{0,1\}^\bN$ with appropriate topologies.

For any $a=\langle a_1,a_2,\ldots\rangle\in \{0,1\}^\bN=T$ consider
the shift map $\sigma$ given by the formula (\cite{RD-intro-to-chaotic}):
$$\sigma(a)=\langle a_2,a_3,\ldots\rangle.$$
Since $T$ is a closed subspace of $[0,1]$, by Tietze's extension theorem
it can be extended to some continuous function $f_0\colon[0,1]\to[0,1]$.

Let
$$S=\left\{\frac{1}{2}\cdot 1, \frac{1}{2}\cdot\frac{1}{3}, \frac{1}{2}\cdot\frac{1}{9}, \ldots\right\};$$
that is, it is the set of centers of most left intervals removed from $[0,1]$ during the classical construction
of the middle-third Cantor set.
Since $\sigma(0)=0$ and $\sigma$ is continuous, we can assume also that
$f_0(s)=0$ for each $s\in S$ (we can multiply original $f_0$ by the continuous function
which is equal to identity on $T$ and equals 0 on $S$).

Let $m=t=1$, $C = [0,1]$, $U=\{0\}$.
Define $f\colon C\times U\to C$ by the formula
$$f(x_0,u_0) = f_0(x_0).$$
Additionally, let
$\zeta^0\in \{0,1\}^\bN = T\subset C$ be given by the formula
$$\zeta^0=\langle 1,0,1,0,0,1,0,0,0,1,0,0,0,0,\ldots\rangle$$
(the sequence of $n$ zeros and one, followed by $n+1$ zeros and one, and so on).
In terms of mapping $(\ref{ex:def-of-homeomorphism})$,
$$\zeta^0 = 2\cdot\sum_{i=2}^\infty \left(\frac{1}{3}\right)^{\frac{i\cdot (i-1)}{2}} \in [0,1].$$
Let $P(x_0)=x_0$ for each $x_0\in [0,1]$,
and $\phi\colon[0,1]\to[0,1]$  be a continuous function
such that $\phi(x_0)=1$ for $x_0\in S\cup\{0\}$, and $\phi(x_0)<1$ otherwise.

Note that for the problem defined in Section~\ref{sec:problem}:
\begin{itemize}
\item $0\in M$ and $S\cap M=\emptyset$;
\item $\zeta^\star=0$, $M^\star=\{\zeta^\star\}$;
\item $D^\star=S\cup\{0\}$.
\end{itemize}
Thus
\begin{enumerate}
\item[$(1)$] the condition $(C1)$ holds: the optimal stationary point $\zeta^\star$ is unique, i.e.~$M^\star=\{\zeta^\star\}$;
\item[$(2)$] the condition $(C3)$ holds: for every $\zeta\in S$ and $u \in U$,
$$P(f(\zeta,u))=f(\zeta,u)=0<\zeta=P(\zeta) 
\textup{\ and\ }
P(f(\zeta^\star,u))=f(\zeta^\star,u)=0=\zeta^\star=P(\zeta^\star).$$
\end{enumerate}

Observe, that for any path $x$ for the system $\probA$:
\begin{itemize}
\item $\Gamma(x)=\{0,\langle 1,0,0,0,\ldots\rangle,\langle 0,1,0,0,0,\ldots\rangle,\langle 0,0,1,0,0,0,\ldots\rangle,\ldots\}$;
  in terms of mapping $(\ref{ex:def-of-homeomorphism})$, $\Gamma(x)=\{0, 2/3, 2/9, 2/27, \ldots \}$.
\item $\Gamma_{\I_d}(x)=\{0\}=\I_d{\rm-}\lim x$.
\end{itemize}

Therefore,
the condition $\condB{\I_d}$ holds (take $x^\star=\langle \zeta^\star,\sigma(\zeta^\star), \sigma(\sigma(\zeta^\star)),\ldots\rangle$),
but $(C2)=\condB{\fin}$ does not hold.  $\Box$
\end{example}

%% file: main.tex
\section{Main results}
\label{sec:main}
The main result of this paper is presented next. The proof of this theorem is provided in Section~\ref{sec:proofs}.

\begin{theorem}\label{thm1:unique}
Suppose that $\I$ is invariant under translations,
$(C1)$, $\condB{\I}$, $(C3)$ hold and $\langle u^{opt},x^{opt} \rangle$
is an optimal process in the problem $\probA, \probB{\I}$.
Then $x^{opt} \to_\I \zeta^\star$,
where $\zeta^\star$ is the unique optimal stationary point from $(C1)$.
\end{theorem}

Note that from part~$(3)$
of Proposition \ref{prop:properties-of-I-convergence} the assertion
``$x^{opt} \to_\I \zeta^\star$'' is equivalent to ``$\Gamma_\I(x)=\{\zeta^\star\}$''.

It is also easy to see that the assertion of Theorem~\ref{thm1:unique}
is true if $D^\star$ is a singleton (i.e.~$D^\star=\{\zeta^\star\}$).
However,
the following example shows that if $\zeta^\star$ is an isolated point of $D^\star$
(if we assume only the first part of condition $(C3)$)
then Theorem~\ref{thm1:unique} may not be true.

\begin{example}\label{ex:dense}
Let $\I=\fin$, $m=1$, $C = U = [0,1]$, and for each $x\in C$:
$$f_0(x)=\left\{\begin{array}{ll}
x & \textup{for\ }  0\leq x\leq\frac{1}{2}+\delta;\\
\frac{\left(0-\left(\frac{1}{2}+\delta\right)\right)\cdot\left(x-\left(\frac{1}{2}+\delta\right)\right)}{1-\left(\frac{1}{2}+\delta\right)}+\left(\frac{1}{2}+\delta\right)  & \textup{for\ }  \frac{1}{2}+\delta<x\leq 1,\\
\end{array}\right.$$
$$f_1(x)=\left\{\begin{array}{ll}
1 & \textup{for\ }  0\leq x\leq\frac{1}{3};\\
\frac{\left(\left(\frac{2}{3}-\delta\right)-1\right)\cdot\left(x-\frac{1}{3}\right)}{\left(\frac{2}{3}-\delta\right)-\frac{1}{3}}+1 & \textup{for\ }  \frac{1}{3}<x\leq\frac{2}{3}-\delta;\\
\frac{\left(\frac{1}{3}-\left(\frac{2}{3}-\delta\right)\right)\cdot\left(x-\left(\frac{2}{3}-\delta\right)\right)}{1-\left(\frac{2}{3}-\delta\right)}+\left(\frac{2}{3}-\delta\right) & \textup{for\ }  \frac{2}{3}-\delta<x\leq 1,\\
\end{array}\right.$$
where $\delta<\frac{1}{12}$ (for example, on Figure~\ref{fig:f}, $\delta=0.05$).
Define $f\colon C\times U\to C$ by the affine formula
$$f(x,u) = f_0(x)\cdot (1-u) + f_1(x)\cdot u.$$
Additionally, let
$\zeta^0=\frac{1}{3}$ and $P(x)=x$ for each $x\in C$.
For the definition of $\phi$ and visualization of $f_0,f_1$ see Figure~\ref{fig:f}.

\tikzstyle{legend_isps}=[rectangle, rounded corners, thin,
                       fill=white!20, text=black, draw]
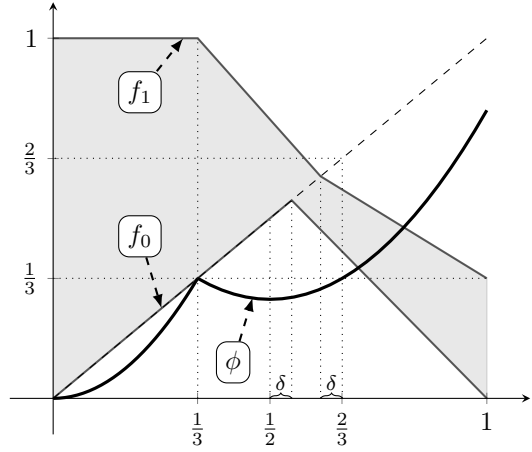
\begin{figure}\centering
\begin{tikzpicture}[scale=1.0]
\begin{axis}[width=8.5cm,xmin=0,xmax=1,xtick={0,0.333,0.5,0.666,1},xticklabels={$0$,$\frac{1}{3}$,$\frac{1}{2}$,$\frac{2}{3}$,$1$},ymin=0,ymax=1,ytick={0,0.333,0.666,1},
yticklabels={$0$,$\frac{1}{3}$,$\frac{2}{3}$,$1$},axis
x line=middle,axis y line=middle,enlargelimits = true]

    \addplot[black,fill=gray!90,opacity=0.2] coordinates {
	(0,0) (1/2+0.05,1/2+0.05) (1,0)
	(1,1/3) (2/3-0.05,2/3-0.05) (1/3,1) (0,1)
    }
    \closedcycle;

    \addplot[black,opacity=0.6,thick]  coordinates {
	(0,0) (1/2+0.05,1/2+0.05) (1,0)
    };
    \addplot[black,opacity=0.6,thick]  coordinates {
	(1,1/3) (2/3-0.05,2/3-0.05) (1/3,1) (0,1)
    };

    \addplot[domain=0:1/3, black, very thick] {x*x*3};
    \addplot[domain=1/3:1, black, very thick] {21/10*x*x-21/10*x+4/5};
    \addplot[domain=0:1, black, dashed] {x};

    \addplot[dotted] plot coordinates {
        (1/3,    0)
        (1/3,    1)
    };
    \addplot[dotted] plot coordinates {
        (0,      1/3)
        (1,      1/3)
    };
    \addplot[dotted] plot coordinates {
        (1/2,    0)
        (1/2,    1/2)
    };
    \addplot[dotted] plot coordinates {
        (2/3,    0)
        (2/3,    2/3)
	(0,      2/3)
    };

    \addplot[dotted] plot coordinates {
        (1/2+0.05,    0)
        (1/2+0.05,    1/2+0.05)
    };
    \addplot[dotted] plot coordinates {
        (2/3-0.05,    0)
        (2/3-0.05,    2/3-0.05)
    };

    \draw [decorate,decoration={brace,amplitude=2pt},color=black] (axis cs: 0.5,0) -- (axis cs: 0.55,0)
       node [midway,anchor=south,inner sep=1pt, outer sep=2pt]{\footnotesize$\delta$};
    \draw [decorate,decoration={brace,amplitude=2pt},color=black] (axis cs: 0.617,0) -- (axis cs: 0.667,0)
       node [midway,anchor=south,inner sep=1pt, outer sep=2pt]{\footnotesize$\delta$};

  \node[legend_isps] (f0) at (axis cs:  0.2,0.45) {$f_0$};
  \node[legend_isps] (f1) at (axis cs:  0.2,0.85) {$f_1$};
  \node[legend_isps] (phi) at (axis cs:  0.415,0.1) {$\phi$};
  \draw[-latex, thick, black, dashed] (f0) -- (axis cs:  0.25,0.25);
  \draw[-latex, thick, black, dashed] (f1) -- (axis cs:  0.3,1);
  \draw[-latex, thick, black, dashed] (phi) -- (axis cs:  0.46,0.28);

\end{axis}
\end{tikzpicture}
\caption{Graph of the quantity $\phi$ and $f_0,f_1$ for Example~\ref{ex:dense}.}\label{fig:f}
\end{figure}

Note that for the problem defined in Section~\ref{sec:problem}:
\begin{itemize}
\item $M=[0,\frac{2}{3}-\delta]$;
\item $\zeta^\star=\zeta^0=\frac{1}{3}$, $\phi^\star=\frac{1}{3}$, $M^\star=\{\frac{1}{3}\}$;
\item $D^\star=\{\frac{1}{3}\}\cup [\frac{2}{3},1]$.
\end{itemize}
Thus
\begin{enumerate}
\item[$(1)$] optimal stationary point $\zeta^\star$ is unique;
\item[$(2)$] the condition $(C2)$ also holds, for example for $x^\star=(\zeta^0,\zeta^0,\zeta^0,\ldots)$, $u^\star=(0,0,0,\ldots)$;
\item[$(3)$] the first part of condition $(C3)$ holds: for every $\zeta\in [\frac{2}{3},1]$ and $u \in U$,
$$P(f(\zeta,u))=f(\zeta,u)\leq f_1(\zeta)<\zeta=P(\zeta).$$
\end{enumerate}
In this example, the process $\langle u^{opt},x^{opt} \rangle,$ where $x^{opt}=(\frac{1}{3},1,\frac{1}{3},1,\frac{1}{3},1,\ldots)$ and
$u^{opt}=(0,1,0,1,0,1,\ldots),$ is an optimal process, however $x^{opt}$ does not converge to $\zeta^\star$ in the sense of $\I$-convergence
(which is equivalent to $\fin$-convergence). $\Box$
\end{example}


Example~\ref{ex:dense} works for classical convergence, statistical convergence,
and for general ideal convergence. It shows that additional assumption
about ``density'' of $D^\star$ in $\zeta^\star$ (i.e.~the second part of condition $(C3)$) is necessary
in \cite{MR3166601}, as well as in \cite{MamPeh-2001-JMAA}.

Recently Leonetti and Caprio in~\cite{MR4428911} proposed another way to bypass the problem indicated in the Example~\ref{ex:dense}:
\begin{itemize}
\item[$(C3{\textrm{-LC}})$:] there exists a linear (and therefore continuous) function $P\colon\bR^m\to\bR$
such that
$$
    P(f(x_0,u_0))<P(x_0)  ~~~~~ {\rm for~all}~~ x_0\in D^\star, u_0\in U, \langle x_0,f(x_0,u_0)\rangle\not=\langle \zeta^\star,\zeta^\star\rangle,
$$
\end{itemize}
where $\zeta^\star$ is an optimal stationary point. It follows from the above condition that $\zeta^\star$ is the unique optimal stationary point,
and it is easy to see that
$(C3{\textrm{-LC}})$ imply $(C3)$. However, we do not have any example of the system with $(C1)+(C3)$ and without $(C3{\textrm{-LC}})$.

\subsection{Special cases}
In this section, we consider two special cases of the ideal convergence; that is, classical convergence and statistical convergence.

\subsubsection*{Classical convergence}
Consider the classical convergence in the problem $\probA$, $\probB{\I}$.
In this case
$$\Gamma(x)=\Gamma_\fin(x)=\left\{a\in\bR^m\colon (x_{n_k})_{k\in\bN}\to a\ \textup{for some subsequence\ }(x_{n_k})\textup{\ of\ }x\right\}$$
is the set of $\omega$-limit points. Condition $\condB{\I}$ is in the form $(C2)$ and functional $\probB{\I}$ is represented in the form
\begin{itemize}
\item[$(**)$:]
$J(x)=J_\fin(x)=\liminf_{k\to\infty} \phi(x_k) \rightarrow \max$.
\end{itemize}

\begin{corollary}\label{cor:3}
Let $(C1)$, $(C2)$, $(C3)$ hold and $\langle u^{opt},x^{opt}\rangle$
is an optimal process in the problem $\probA, (**)=\probB{\fin}$.
Then $x^{opt}$ converges to $\zeta^\star.$
\end{corollary}

\subsubsection*{Statistical convergence}
Now consider the statistical convergence instead of ideal convergence in the problem $\probA, \probB{\I}$.
Functional $\probB{\I}=\probB{\I_d}$ in this case can be defined as follows
\begin{itemize}
\item[$\probB{\I_d}$:]
$J_{\I_d}(x)=\C-\liminf_{k\to\infty}\phi(x_k) \rightarrow \max,$
\end{itemize}
where
$\C-\liminf_{k\to\infty}\phi(x_k)=\I_d-\liminf x$ stands for the minimal element in the set of statistical cluster points.
Recall also that according to Example \ref{ex:statistical-not-fin}, condition $(C2)$ is stronger than $\condB{\I_d}$.

\begin{corollary}\label{cor:4}
Let $(C1)$, $\condB{\I_d}$, $(C3)$ hold and $\langle u^{opt},x^{opt}\rangle$
is an optimal process in the problem $\probA, \probB{\I_d}$.
Then $x^{opt}$ statistically converges to $\zeta^\star.$
\end{corollary}

%% file: proofs.tex
\section{Proof of Theorem \ref{thm1:unique}}
\label{sec:proofs}


For every $r\in\bR$
define the set
$$D_r=\left\{\zeta\in C\colon ~\phi(\zeta)\geq r\right\}.$$
Clearly, $D^\star=D_{\phi^\star}$.
For any continuous function $P\colon\bR^m\to\bR$ let
$$E_P=\left\{\zeta\in\bR^m\colon P(f(\zeta,u_0))<P(\zeta)\textup{\ for all\ }u_0\in U\right\},$$
and
$$\overline{E}_P=\left\{\zeta\in\bR^m\colon P(f(\zeta,u_0))\leq P(\zeta)\textup{\ for all\ }u_0\in U\right\}.$$
It is clear that $M\cap E_P=\emptyset$.
If $A\subset\bR^m$ is compact then
$$\arg\min_{\zeta\in A}P(\zeta) \dot{=} \left\{\zeta_1\in A\colon P(\zeta_1)=\min_{\zeta\in A} P(\zeta)\right\}.$$
Analogously we define operator $\arg\max$.

\begin{lemma}\label{le:claim1}
	Assume that $\I$ is invariant under translations,
        $r\in\bR$ and $\langle u,x\rangle$ is a process in the problem $\probA, \probB{\I}$
	with $J_\I(x)\geq r$. If $P\colon\bR^m\to\bR$ is a continuous function
	then
	$$\arg\min_{\zeta\in\Gamma_\I(x)}P(\zeta)\subset D_r\setminus E_P.$$
\end{lemma}

{\bf Proof.~}
As $J_\I(x)\geq r$,
by Lemma~\ref{le:functional-j},
$J_\I(x)=\min_{\zeta\in\Gamma_\I(x)}\phi(\zeta)\geq r$.
Thus $\Gamma_\I(x)\subset D_r$, and so $\arg\min_{\zeta\in\Gamma_\I(x)}P(\zeta)\subset D_r$.

Let
$$F(\zeta)=\max_{u_0\in U} P(f(\zeta,u_0))-P(\zeta).$$

It is clear that $F\colon \bR^m\to\bR$ is continuous, and
\begin{equation}\label{claim1:02}
   F(\zeta)<0 ~~~~ {\rm for ~ all }~~ \zeta\in E_P.
\end{equation}


Suppose that there exists $\zeta_1\in\Gamma_\I(x)$
such that $\zeta_1\in E_P$ and $$\min_{\zeta\in\Gamma_\I(x)} P(\zeta)=P(\zeta_1).$$
Denote $\delta = - F(\zeta_1)/8.$ Clearly $\delta>0$ thanks to (\ref{claim1:02}).

Since functions $F$ and $P$ are continuous and $\Gamma_\I(x)$ is a compact set, there exists $\gamma>0$ such that
\begin{equation}\label{claim1:align1}
    \forall_{\zeta\in B(\zeta_1,\gamma)} F(\zeta) \leq -4\delta, \textup{\ and }
\end{equation}
\begin{equation}\label{claim1:align2}
    \forall_{\zeta\in B(\zeta_1,\gamma)} P(\zeta) \leq P(\zeta_1)+\delta, \textup{\ and }
\end{equation}
\begin{equation}\label{claim1:align3}
    \forall_{\zeta\in B(\Gamma_\I(x),\gamma)} P(\zeta) \geq \min_{y\in \Gamma_\I(x)}P(y)-\delta.
\end{equation}

If $x_k\in B(\zeta_1,\gamma)$ then $F(x_k)\leq-4\delta$, i.e.
$P(f(x_k,u_0))\leq P(x_k)-4\delta$
for each $u_0\in U$ and in particular for $u_k\in U$ that leads to
$P(x_{k+1})\leq P(x_k)-4\delta.$ Moreover, from (\ref{claim1:align2}) we have
$P(x_k)\leq P(\zeta_1)+\delta$ and therefore
$$P(x_{k+1}) \leq P(\zeta_1)-3\delta.$$

On the other hand (\ref{claim1:align3}) implies
$$\forall_{\zeta\in B(\Gamma_\I(x),\gamma)} P(\zeta)\geq \min_{y\in \Gamma_\I(x)}P(y)-\delta
= P(\zeta_1) -\delta > P(\zeta_1)-3\delta.$$
Thus $x_{k+1}\notin B(\Gamma_\I(x),\gamma)$.
By the above considerations we get
\begin{equation}\label{claim1:eq:implication}
x_k\in B(\zeta_1,\gamma) \Longrightarrow x_{k+1}\notin B(\Gamma_\I(x),\gamma).
\end{equation}
This contradicts with Lemma~\ref{le:cluster-ideal-invariant}.
$\Box$ \bigskip


\begin{lemma}\label{le:claim2}
	Assume that $\I$ is invariant under translations,
	$r\in\bR$ and $\langle u,x\rangle$ is a process in the problem $\probA, \probB{\I}$
	with $J_\I(x)\geq r$. If $P\colon\bR^m\to\bR$ is a continuous function
	and $D_r\setminus E_P\subset\overline{E}_P$
	then
	$$\arg\max_{\zeta\in\Gamma_\I(x)}P(\zeta)\cap\left( D_r\setminus E_P\right) \not=\emptyset.$$
\end{lemma}

{\bf Proof.~}
As in the proof of Lemma~\ref{le:claim1} observe that
$\Gamma_\I(x)\subset D_r$, $\arg\max_{\zeta\in\Gamma_\I(x)}P(\zeta)\subset D_r$, and define
$$F(\zeta)=\max_{u_0\in U} P(f(\zeta,u_0))-P(\zeta).$$
Again, $F\colon \bR^m\to\bR$ is continuous, and
\begin{equation}\label{claim2:02:01}
    F(\zeta)<0  ~~{\rm for ~ all }~~ \zeta\in E_P;
\end{equation}
\begin{equation}\label{claim2:02:02}
    F(\zeta)\geq 0  ~~{\rm for ~ all }~~ \zeta \notin E_P;
\end{equation}
\begin{equation}\label{claim2:02:03}
    F(\zeta)= 0  ~~{\rm for ~ all }~~ \zeta \in \overline{E}_P\setminus E_P.
\end{equation}
The last equality follows from the previous ones and the fact that $F$ is continuous.

Denote
$$Z_1 =   \Gamma_\I(x)\setminus E_P\subset D_r\setminus E_P, ~~
Z_2 = \arg\max_{\zeta\in\Gamma_\I(x)}P(\zeta),$$
%


\noindent
and assume (contrary to the lemma assertion) that $Z_1\cap Z_2=\emptyset$.
Note that $Z_1,Z_2$ are compact and,
from the assumption that $Z_1\cap Z_2=\emptyset$ it follows that
$\max_{\zeta\in Z_1}P(\zeta)<\min_{\zeta\in Z_2}P(\zeta)$
(in fact, $P\restriction Z_2$ is constant and equal
to the maximum value of $P$ on $\Gamma_\I(x)$; if it is equal to $\max_{\zeta\in Z_1}P(\zeta)$
then it follows from the definition of $Z_2$ that $Z_1\cap Z_2\not=\emptyset$).
Since (by the assumption of the lemma) $D_r\setminus E_P\subset\overline{E}_P$, $(\ref{claim2:02:03})$ gives us $F\restriction Z_1=0$.

Denote also
$$p_1 = \max_{\zeta\in Z_1}P(\zeta), ~~ p_2 = \min_{\zeta\in Z_2}P(\zeta).$$
%
Let
\begin{equation}\label{03}
    a = \frac{p_2-p_1}{8} > 0.
\end{equation}

{\bf 1.~} Since functions $F,$ $P$ are continuous and $F\restriction Z_1=0$,
there exists $\gamma>0$ such that
\begin{equation}\label{align11}
    \forall_{\zeta\in B(Z_1,\gamma)} F(\zeta) \leq 4 a, \textup{\ and }
\end{equation}
\begin{equation}\label{align12}
    \forall_{\zeta\in B(Z_1,\gamma)} P(\zeta) \leq p_1+a, \textup{\ and }
\end{equation}
\begin{equation}\label{align13}
    \forall_{\zeta\in B(Z_2,\gamma)} P(\zeta) \geq p_2-a.
\end{equation}

Let $x_{k-1} \in B(Z_1,\gamma)$. Then from (\ref{align11}), (\ref{align12}) the following two relations hold

$$P(x_k) - P(x_{k-1}) = P(f(x_{k-1}, u_{k-1})) - P(x_{k-1}) \le F(x_{k-1}) \le 4 a;$$
$$P(x_{k-1}) \le p_1 + a.$$
From these inequalities we have
$$P(x_k) \le p_1 + 5 a .$$
From (\ref{03}) it follows $p_1 =  p_2 - 8 a$ and then
$$P(x_k) \le p_2 - 8 a + 5 a < p_2 - a.$$
According to (\ref{align13}) this means that $x_k \notin B(Z_2,\gamma)$.
Therefore we conclude that
\begin{equation}\label{eq:implication2}
x_k\in B(Z_2,\gamma) \Longrightarrow x_{k-1}\notin B(Z_1,\gamma).
\end{equation}

\bigskip

{\bf 2.~} We fix the number $\gamma$ and consider the set
\begin{equation}\label{04}
   \Gamma ~ = ~  \Gamma_\I(x) \setminus B(Z_1,\gamma).
\end{equation}
From~(\ref{03}), (\ref{align12}) and~(\ref{align13}) and the fact that $Z_2\subset\Gamma_\I(x)$
it follows that $Z_2\subset\Gamma$.
Moreover, $\Gamma\subset D_r\cap E_P$,
and~$(\ref{claim2:02:01})$ implies that $F(\zeta) < 0$ for all $\zeta \in \Gamma$.
Denote
$$\delta = - \max_{\zeta \in \Gamma} F(\zeta) > 0.$$
Take any number $\varepsilon>0$ satisfying 
\begin{equation}\label{eq:4-epsilon-delta}
4 \varepsilon < \delta.
\end{equation}

Since functions $F,$ $P$ are continuous
there exists a sufficiently small number $\eta \in (0, \gamma)$ such that
\begin{equation}\label{align21}
    \forall_{\zeta\in B(\Gamma,\eta)} F(\zeta) \leq - \delta + \varepsilon,  \textup{\ and }
\end{equation}
\begin{equation}\label{align22}
    \forall_{\zeta\in B(Z_2,\eta)} P(\zeta) \geq p_2 - \varepsilon, \textup{\ and }
\end{equation}
\begin{equation}\label{align23}
    \forall_{\zeta\in B(\Gamma_\I(x),\eta)} P(\zeta) \leq p_2 + \varepsilon.
\end{equation}

We show by contradiction that there is no $k$ such that 
\begin{equation}\label{k-contradiction}
x_k\in B(Z_2,\eta) \textup{\ and\ } x_{k-1}\in B(\Gamma,\eta).
\end{equation}
Suppose that $k$ fulfills $(\ref{k-contradiction})$.
From (\ref{align21}) we have
$$P(x_k) - P(x_{k-1}) \le F(x_{k-1}) \le - \delta + \varepsilon$$
or
$$P(x_{k-1}) \ge P(x_k) + \delta - \varepsilon.$$
Then, from (\ref{align22}) it follows
$$P(x_{k-1}) \ge (p_2 - \varepsilon) + \delta -  \varepsilon =p_2 + \delta - 2\varepsilon,$$
and by (\ref{eq:4-epsilon-delta})
$$P(x_{k-1}) > p_2 +  2 \varepsilon.$$

On the other hand (\ref{align23}) yields
$$P(x_{k-1}) \le p_2 + \varepsilon.$$

The last two inequalities lead to a contradiction. This proves that
the relations $x_k\in B(Z_2,\eta)$ and $x_{k-1}\in B(\Gamma,\eta)$ cannot be satisfied at the same time. Therefore
the following is true:
\begin{equation}\label{eq:implication3}
x_k\in B(Z_2,\eta) \Longrightarrow x_{k-1}\notin B(\Gamma,\eta).
\end{equation}

\bigskip

{\bf 3.~} Now since $\eta < \gamma$, it is not difficult to observe that
the relation
$$B(\Gamma_\I(x),\eta) \subset B(\Gamma,\eta) \cup B(Z_1,\gamma)$$
holds. Then (\ref{eq:implication2}) and (\ref{eq:implication3}) implies that
\begin{equation}\label{eq:implication4}
x_k\in B(Z_2,\eta) \Longrightarrow x_{k-1}\notin B(\Gamma_\I(x),\eta).
\end{equation}
The above implication contradicts with Lemma~\ref{le:cluster-ideal-invariant}.
$\Box$ \bigskip


{\bf Proof of Theorem \ref{thm1:unique}.~~}
Let $r=\phi^\star$. Then $D_r=D_{\phi^\star}=D^\star$.
By $\condB{\I}$ for the process $\langle u^\star,x^\star\rangle$, $J_\I(x^\star)\geq r$.

Fix the function $P$ like in $(C3)$. Then 
$$D_r\setminus E_P=D^\star\setminus E_P\subset M^\star\subset \overline{E}_P\setminus E_P.$$

Since $J_\I(x^\star)=r$, the maximal
value of the functional $\probB{\I}$ is not less than $r$.
As $\langle u^{opt}, x^{opt} \rangle$ is an optimal process,
$J_\I(x^{opt})\geq r$.
Thus, by Lemma~\ref{le:claim1},
$$\arg\min_{\zeta\in\Gamma_\I(x^{opt})} P(\zeta) \subset D^\star\setminus E_P \subset M^\star=\{\zeta^\star\},$$
where $\zeta^\star$ is the unique optimal stationary point from $(C1)$.
Then, by Lemma~\ref{le:claim2},
$$\zeta^\star\in \arg\max_{\zeta\in\Gamma_\I(x^{opt})} P(\zeta) \cap M^\star.$$
Thus $P(\zeta)=P(\zeta^\star)$ for all $\zeta\in\Gamma_\I(x^{opt})$.
It follows that
$$\Gamma_\I(x^{opt})=\arg\min_{\zeta\in\Gamma_\I(x^{opt})} P(\zeta)\subset M^\star=\{\zeta^\star\}.$$
From part~$(3)$
of Proposition \ref{prop:properties-of-I-convergence} we obtain $x^{opt}\to_{\I}\zeta^\star$.
$\Box$ \bigskip